\documentclass{amsart2000}
\theoremstyle{definition}
\newtheorem{defin}{Definition}[section]
\newtheorem{example}[defin]{Example}
\newtheorem{remark}[defin]{Remark}
\theoremstyle{plain}
\newtheorem{lemma}[defin]{Lemma}
\newtheorem{prop}[defin]{Proposition}
\newtheorem{theorem}[defin]{Theorem}
\newtheorem{corol}[defin]{Corollary}

\newcommand{\ra}{\rightarrow}
\newcommand{\la}{\leftarrow}
\newcommand{\ua}{\uparrow\!}

\newcommand{\se}{\subseteq}
\newcommand{\T}{\mathcal{T}}
\newcommand{\R}{\mathcal{R}}
\newcommand{\Ra}{\Rightarrow}
\newcommand{\LRa}{\Leftrightarrow}
\newcommand{\adj}{^{\dashv}}
\newcommand{\Q}{\mathcal{Q}}
\newcommand{\Two}{\mathbf{2}}
\renewcommand{\P}{\mathcal{P}}
\newcommand{\I}{\bigwedge}
\newcommand{\lb}{\langle}
\newcommand{\rb}{\rangle}
\newcommand{\B}{\mathcal{B}}
\newcommand{\x}{\cdot}
\renewcommand{\S}{\bigvee}
\renewcommand{\_}{\text{\textvisiblespace\hspace{.07em}}}
\newcommand{\op}{^{\mathrm{op}}}
\DeclareMathOperator{\Max}{Max}

\begin{document}
\noindent
Topology Atlas Invited Contributions \textbf{7} no. 2 (2002) 13 pp.
\vspace{0.25in}
\title{On simple and semisimple quantales}
\author{David Kruml}
\address{Department of Mathematics, Masaryk University\\
Jan\'{a}\v{c}kovo n\'{a}m. 2a, 635 00 Brno}
\email{kruml@math.muni.cz}
\author{Jan Paseka}
\email{paseka@math.muni.cz}
\thanks{Financial Support of the Grant Agency of the Czech Republic under 
the grant No.~201/99/0310 is gratefully acknowledged.}
\thanks{The authors thank Mamuka Jibladze for suggestive questions given 
during Quantum logic/quantale/category workshop in Brussels 10--14/4/2000, 
from which the paper arose.}
\thanks{This article was accepted for publication in \textit{Proceedings 
of the Ninth Prague Topological Symposium, (Prague, 2001)}, 
edited by Petr Simon, but was omitted from the final 
publication due to a production error.}
\begin{abstract}
In a recent paper \cite{PeRo97}, J.~W.~Pelletier and J.~Rosick\'{y} 
published a characterization of *-simple *-quantales. 
Their results were adapted for the case of simple quantales by J.~Paseka 
in \cite{Pa97}. 
In this paper we present similar characterizations which do not use a 
notion of discrete quantale. 
We also show a completely new characterization based on separating and 
cyclic sets. Further we explain a link to simple quantale modules. 
To  apply these characterizations, we study (*-)semisimple (*-)quantales 
and discuss some other perspectives. 
Our approach has connections with several earlier works on the subject 
\cite{Pa97,PeRo97}.
\end{abstract}
\maketitle

\section{Preliminaries}

\begin{defin}
By a {\em sup-lattice} is meant a complete lattice, a {\em sup-lattice
morphism} is a mapping preserving arbitrary joins.

If $f\colon S\ra T$ is a sup-lattice morphism, the assignment 
$$f(s) \leq t\LRa s \leq f\adj(t),$$ 
explicitly
$$f\adj(t) = \S\{ s \mid f(s) \leq t\},$$ 
defines a mapping $f\adj \colon T\ra S$ preserving all meets. 
This mapping  $f\adj$ is called the {\em adjoint} of $f$.

The top element of a sup-lattice is denoted by $1$, the bottom element by 
$0$.
\end{defin}

\begin{defin} 
By a {\em quantale} is meant a sup-lattice $Q$ equipped with an 
associative multiplication which distributes over joins
\begin{align*}
a\left(\S a_i\right)&=\S(aa_i),&
\left(\S a_i\right)a&=\S(a_ia)
\end{align*}
for all $a,a_i\in Q$.

An element $a\in Q$ is called {\em right-sided} or {\em left-sided} if 
$a1 \leq a$ or $1a \leq a$, respectively. 
We write $a\in\R(Q)$, resp. $a\in\mathcal{L}(Q)$.
Elements of $\T(Q)=\R(Q)\cap\mathcal{L}(Q)$ are called {\em two-sided}.

A {\em quantale morphism} is a sup-lattice morphism preserving 
multiplication.

An equivalence relation $\sim$ on $Q$ is said to be a {\em congruence} if
\begin{align*}
a\sim b&\Ra ca\sim cb,ac\sim bc,&
a_i\sim b_i&\Ra\S a_i\sim\S b_i
\end{align*}
for all $a,b,c,a_i,b_i\in Q$. 
One can easily check that congruences are exactly coset equivalences given 
by quantale morphisms.

Actions $a\_$, $\_a$ for fixed $a\in Q$ determine sup-lattice 
endomorphisms on $Q$. Their adjoints are denoted by $\_\la a$, $a\ra\_$, 
respectively.\footnote{The notation of residuations in literature varies. 
Rosenthal \cite{Rose90} uses $\ra_l,\ra_r$, our notation is motivated by 
the first assertion of \ref{resid}, which makes possible to omit 
parenthesis.} 
That is, 
$${a\leq b\ra c} \LRa {ab\leq c} \LRa {b\leq c\la a}.$$
The operations $\ra$, $\la$ are called {\em residuations}.

A set $P\se Q$ is said to be {\em residually closed} if $a\ra p\in P,p\la 
a\in P$ for every $a\in Q,p\in P$.
\end{defin}

\begin{lemma}\cite{Rose90} 
Let $Q$ be a quantale, $a,b,c,a_i\in Q$. 
Then
\begin{align*}
(b\ra c)\la a&= b\ra (c\la a),\\
a\ra (b\ra c)&= (ab)\ra c,\\
(c\la a)\la b&= c\la (ab),\\
\left(\S a_i\right)\ra c&= \I(a_i\ra c),\\
c\la\left(\S a_i\right)&= \I(c\la a_i).
\end{align*}\label{resid}
\end{lemma}

\begin{defin}
By a {\em *-quantale}\footnote{In \cite{MuPe00,PeRo97,Kr00}, *-quantales 
are called involutive quantales, as well as the following concepts.}
is meant a quantale $Q$ with a unary operation of {\em involution} such that
\begin{align*}
a^{**}&= a,&
(ab)^*&= b^*a^*,&
\left(\S a_i\right)^*&= \S a_i^*
\end{align*}
for all $a,b,a_i\in Q$.

An element $a\in Q$ is called {\em hermitian} if $a^*=a$, the set of 
hermitian elements is denoted by $\mathcal{H}(Q)$.
 
By a {\em *-quantale morphism} (or simply a {\em *-morphism}) is meant a 
quantale morphism of *-quantales which also preserves the involution.

A quantale congruence $\sim$ on a *-quantale $Q$ is said to be a {\em 
*-congruence} if $$ a\sim b\Ra a^*\sim b^* $$ for all $a,b\in Q$.
\end{defin}

Note that $(a\ra b)^*=b^*\la a^*$.

\begin{defin}
Let $Q$ be a quantale. 
By a left (right) $Q$-module is meant a sup-lattice $M$ with an action 
$Q\times M\to M$ ($M\times Q\to M$) provided that 
\begin{align*}
\begin{aligned}
a\left( \S m_i\right)&= \S (am_i),\\
\left( \S a_i\right) m&= \S (a_im),\\ 
a(bm)&= (ab)m 
\end{aligned}
& &
\left(\begin{aligned}
\left( \S m_i\right) a&= \S (m_ia),\\
m\left(\S a_i\right)&= \S (ma_i),\\
(mb)a &= m(ba)
\end{aligned}\right)
\end{align*}
for all $a,b,a_i\in Q, m,m_i\in M$.

An adjoint of the module action is called a residuation as well and 
denoted by $\la$ ($\ra$). 
One can easily verify similar properties of the residuation as in 
\ref{resid}, from whose it follows that $M\op$ with the residuation is a 
right (left) $Q$-module. 
$M\op$ is called a {\em dual module} of $M$.\footnote{The dual module, as  
well as submodules and quotient modules, is determined uniquely up to 
isomorphism. 
For simplicity, we ignore the difference of isomorphic modules.}

A mapping $f\colon M\to N$ of two left (right) $Q$-modules $M,N$ is said 
to be 
a {\em $Q$-module morphism} if $f$ is a sup-lattice morphism and 
$f(am)=af(m)$ ($f(ma)=f(m)a$) for every $a\in Q,m\in M$. Let us recall 
that the adjoint $f\adj \colon N\op\to M\op$ of a $Q$-module morphism 
$f\colon M\to N$ 
is a $Q$-module morphism of dual modules (cf. \cite{AbVi93}). Two special 
cases of the module morphism is a submodule inclusion and a quotient 
mapping. An adjoint of
inclusion is a quotient mapping and vice versa.

A $Q$-module $M$ is called {\em faithful} if
$(\forall m\in M)(am=bm)\Ra a=b$ for all $a,b\in Q$.
\end{defin}

Through the paper we do not distinguish between left and right modules if 
it is not important or if the meaning is clear. 
We write $1_Q$ and $1_M$ for the quantale and module top, respectively, 
whenever a confusion is imminent.

\section{Simple and semisimple quantales}

The original construction \cite{PeRo97,Pa97} characterizes simple
quantales as discrete faithful factors. We are using stronger definition
of faithfulness, in which discreteness is enclosed. In fact, this idea
appeared in the proof of Pelletier and Rosick\'y but it was not fully
exploited.

\begin{defin}
By a {\em simple quantale} is meant a quantale $Q$ such that 
$1\cdot
1\neq 0$ and every its quotient morphism is either an isomorphism or a constant
morphism on a trivial quantale. In other words, there are exactly two
congruences on $Q$.

By a {\em semisimple quantale} is meant a subdirect product of simple
quantales. That is, a quantale is semisimple if it has enough surjective
morphisms onto simple quantales to separate elements.
\end{defin}

\begin{defin}
A quantale $Q$ is called a {\em factor} if $\T (Q)=\{ 
0,1\}$ and $0\neq 1$.

$Q$ is said to be {\em strictly faithful} if
$(\forall l\in\mathcal{L}(Q),\forall r\in\R (Q))(lar=lbr)\Ra a=b$
for all $a,b\in Q$.
\end{defin}

\begin{defin}
Let $Q$ be a quantale.
A set $P\se Q$ is said to be {\em separating} if for every 
$a,b\in Q,a\not\geq b$ there is an element $p\in P$ provided that
$a\leq p,b\not\leq p$. Note that then every element of $Q$ is a meet of
elements from $P$.

A nonempty set $P,1\not\in P$ is called {\em cyclic} if $\{ a\ra p\la b\mid
a,b\in Q\}=P\cup\{ 1\}$ for every $p\in P$. In fact, any
quantale having a cyclic set is nontrivial.
\end{defin}

\begin{lemma}\label{cyclic}
Let $Q$ be a quantale.
A nonempty set $P\se Q-\{1\}$ is cyclic, iff $P\cup\{ 1\}$ is residually 
closed, $t=\I P=1\ra p\la 1$ for every $p\in P$, and every element 
$p\in P$ can be recovered from $t$, i.e., $p=a\ra t\la b$ for some $a,b\in Q$.
Then $t$ is two-sided and $a,b$ can be taken right-sided or left-sided,
respectively.
\end{lemma}

\begin{proof}
Let $P$ be cyclic and $p\in P$. 
Then
$$
1\not =
t = 
\I P =
\I_{a,b\in Q}(a\ra p\la b) = 
\biggl(\S_{\ a\in Q}a\biggr)\ra p\la\biggl(\S_{\ b\in Q}b\biggr) = 
1\ra p\la 1,
$$
i.e., $t\in P$. 
Thus $P\cup\{ 1\}=\{ a\ra t\la b\mid a,b\in Q\}$.
Since $1\ra t\la 1\in P\cup\{ 1\}$ and $1\ra t\la 1\leq a\ra t\la b$ for 
every $a,b\in Q$, we have $t=1\ra t\la 1$. 
Then 
\begin{equation*}
a\ra t\la b
= a\ra(1\ra t\la 1)\la b
= (a1)\ra t\la(b1)
\end{equation*}
and 
$a1\in\R(Q)$,$1b\in\mathcal{L}(Q)$. 
Finally, 
\begin{equation*}
\begin{split}
t&
= 1\ra t\la 1\\
&
\leq (1\x 1)\ra t\la 1\\
&
= 1\ra(1\ra t\la 1)\\
&
= 1\ra t, 
\end{split}
\end{equation*} 
hence $t1\leq t$. 
Similarly, $1t\leq t$. 
It is evident that $P\cup\{ 1\}$ is residually closed.

Conversely, let $P$ satisfies the considered properties. 
Then for every $p\in P$ we have $p=a\ra t\la b$ for some $a,b\in Q$.
Thus 
\begin{align*}
P'&
= \{ c\ra p\la d\mid c,d\in Q\}\\
& 
= \{(ca)\ra t\la(bd) \mid c,d\in Q\}\\
&
\se P\cup\{ 1\}\\
\intertext{and}
P\cup\{ 1\}&
=\{ c\ra t\la d\mid c,d\in Q\}\\
& 
= \{(c1)\ra p\la(1d)\mid c,d\in Q\}\\
&
\se P',
\end{align*} 
i.e., $P$ is cyclic.
\end{proof}

\begin{theorem}\label{simple}
Let $Q$ be a quantale. 
The following conditions are equivalent:
\begin{enumerate}
\item
$Q$ is a simple quantale.
\item
$Q$ is a strictly faithful factor.
\item
$Q$ has a separating cyclic set.
\end{enumerate}
\end{theorem}

\begin{proof}
(1)$\Ra$(2): 
Let $t\in\T (Q)$. 
The assignment $a\mapsto a\vee t$ defines a surjective morphism 
$Q\ra\ua t$. 
Since $Q$ is simple, we have either $t=0$ or $t=1$, i. e. $Q$ is a factor. 
From the additional condition $1\x 1\neq 0$ we conclude $1\x 1=1$, because 
$1\x 1$ is always two-sided. 
Let us put
$$a\sim b\LRa(\forall l\in\mathcal{L}(Q),\forall r\in\R 
(Q))(lar=lbr).$$
The relation $\sim$ is a quantale congruence. 
Indeed, $cr\in\R(Q)$ and $lc\in\mathcal{L}(Q)$ for every 
$c\in Q,r\in\R(Q),l\in\mathcal{L}(Q)$ and hence $a\sim b\Ra ca\sim cb$ 
and 
$ac\sim bc$. 
It is not hard to see that $\sim$ is a sup-lattice congruence. 
If $0\sim 1$, then $0=1\x 0\x 1=1\x 1\x 1=1$. 
Thus $\sim$ is necessarily a diagonal relation, i.e., $a\sim b\LRa a=b$, 
and so $Q$ is strictly faithful.

(2)$\Ra$(3): 
Put $P=\{ r\ra 0\la l\mid r\in\R (Q),l\in\mathcal{L}(Q)\} -\{ 1\}$. 
$P\cup\{ 1\}$ 
is residually closed, because $b\ra (r\ra 0\la l)\la a=(br)\ra 0\la (la)$. 
From the strict faithfulness we have that $1\x 1\x 1\neq 1\x 0\x 1$. 
Since $Q$ is factor and $1\x 1\in\T(Q)$, we have $1\x 1=1$. 
For $p\in P$ we have $1\ra p\la 1\in\T(Q)=\{ 0,1\}$. If $1\ra p\la 1=1$, 
then $1=1\x 1\x 1\leq p$, thus $1\ra p\la 1=0$ and hence $P$ is cyclic. 
Let $a\not\geq b$, then there exist $r\in\R (Q),l\in\mathcal{L}(Q)$ 
such that 
$lar\not\geq lbr$. 
Since $Q$ is a factor and $lar,lbr$ are two-sided, $lar=0$ and $lbr=1$. 
But it means that $a\leq r\ra 0\la l$ and $b\not\leq r\ra 0\la l$, i.e., 
$P$ is separating.

(3)$\Ra$(1): 
Let $P$ be separating and cyclic. 
Then $0=\I P=1\ra p\la 1$ for every $p\in P$ and there are $a,b\in Q$ such 
that $p=b\ra 0\la a$. 
Let $a\in Q$ and suppose $1a1\leq p$ for some $p\in P$. 
Then $a\leq 1\ra p\la 1=0$, that is $a\neq 0\Ra 1a1=1$. 
Further $0\not=1=1\x 1\x 1\leq 1\x 1$.
Consider $a\not\geq b$ and a congruence $\sim$ generated by $a\sim b$. 
There exists a $p\in P,p=d\ra 0\la c$ such that $a\leq p,b\not\leq p$. 
Then $cad=0,cbd\neq 0,cad\sim cbd$ and thus $0=1cad1\sim 1cbd1=1$, i.e.,
there are only two congruences on $Q$, so that $Q$ is simple.
\end{proof}

Let us recall from \cite{Kr00} the following definition.

\begin{defin}
An element $p$ of a quantale $Q$ is said to be {\em prime} if
$$a1b\leq p\Ra a\leq p
\text{ or }
b\leq p $$
for all $a,b\in Q$. 
Let $\P(Q)$ denote the set of all primes of $Q$.
\end{defin}

\begin{theorem}
Let $Q$ be a quantale. Then there is a bijective correspondence between 
cyclic sets and two-sided primes.
\end{theorem}

\begin{proof}
Let $p$ be a two-sided prime. 
We show that $P=\{ a\ra p\la b\mid a,b\in Q\}$ is cyclic. 
It is sufficient to prove that $1\ra(a\ra p\la b)\la 1=p$ if 
$a\ra p\la b\neq 1$, because 
$c\ra(a\ra p\la b)\la d=(ca)\ra p\la(bd)\in P$ for every $c,d\in Q$. 
Note that from the two-sidedness of $p$ it follows that 
$a1b\leq p \LRa a\leq p$ or $b\leq p$. 
Then 
\begin{equation*}
\begin{split}
q\leq 1\ra(a\ra p\la b)\la 1&
\LRa b1q1a\leq p\\
&
\LRa b1q\leq p \text{ or } a\leq p\\
&
\LRa b\leq p \text{ or } q\leq p \text{ or } a\leq p.
\end{split}
\end{equation*}
But from $b\leq p$ or $a\leq p$ we have $b1a\leq p$, i.e., 
$a\ra p\la b=1$. Thus
$q\leq 1\ra(a\ra p\la b)\la 1\LRa q\leq p$, 
i.e., $1\ra(a\ra p\la b)\la 1=p$.

Conversely, let $P$ be a cyclic set, $p\in P$. 
Put $t=1\ra p\la 1$, then evidently $t\in\T(Q)$. 
Since $P$ is cyclic and  $1(1\ra t)1\leq 1t\leq t$, we have that 
$t\leq 1\ra t\leq 1\ra t\la 1=t$, i.e., $1\ra t=t$ and similarly  
$t\la 1=t$. 
Suppose now that $a1b\leq t$ and $a\not\leq t$. 
Then $a\not\leq 1\ra t$, i.e., $t\la a\neq 1$ and thus 
\begin{equation*}
\begin{split}
t&
=1\ra (t\la a)\la 1\\
&
=(1\ra t\la a)\la  1\\
&
=(t\la a)\la 1. 
\end{split}
\end{equation*}
Hence 
\begin{equation*}  
\begin{split}
b\ra t&
=b\ra(t\la a)\la 1\\
&
=(b\ra t\la a)\la 1\\
&
=1\la 1\\
&
=1,
\end{split}
\end{equation*}
i.e., $b\leq t\la 1=t$. We have proved that $t$ is prime.
\end{proof}

\begin{lemma}\label{primes_are_closed}
If $p\in\P(Q)$ and $a\in Q$, then $a\ra p,p\la a\in\P(Q)\cup\{1\}$.
\end{lemma}

\begin{proof}
\begin{equation*}
\begin{split}
b1c\leq a\ra p&
\Ra b1\cdot 1c\leq a\ra p\\
&
\LRa b1\cdot 1ca\leq p\\
&
\Ra b1\leq p \text{ or } ca\leq p\\
&
\Ra b\leq 1\ra p\leq a\ra p \text{ or } c\leq a\ra p.
\end{split} 
\end{equation*}  
Similarly for $p\la a$.
\end{proof}

\begin{corol}
Elements of cyclic sets are primes.
\end{corol}

\begin{prop}\label{simple-quotients}
Let $Q$ be a quantale, $P$ a cyclic set. For $a,b\in Q$ put
$$a\sim b\LRa (\forall p\in P)(a\leq p\LRa b\leq p).$$
Then $\sim$ is a quantale congruence with the property that the quotient
quantale $Q/\sim$ is simple. Conversely, if $f\colon Q\ra K$ is a 
surjective 
morphism onto simple quantale $K$ with separating cyclic set $R$, then 
$f\adj[R]$ is a cyclic set in $Q$.
\end{prop}

\begin{proof}
$\sim$ is evidently a sup-lattice congruence.
Suppose that $a\sim b$, i.e., $a\leq p\LRa b\leq p$ for every $p\in P$. 
But then
also $a\leq c\ra p\LRa b\leq c\ra p$ and $a\leq p\la c\LRa b\leq p\la c$, i.e.
$ac\sim bc$ and $ca\sim cb$.

Let $g\colon Q\ra Q/\sim$ be the corresponding quotient morphism. 
By definition, $g[P]$ is separating and $g\adj g|_P=id_P$. 
Then 
\begin{equation*}
\begin{split}
g(b)\leq g(a)\ra g(p)&
\LRa g(ba)=g(b)g(a)\leq g(p)\\
&
\LRa ba\leq g\adj g(p)=p\\
&
\LRa b\leq a\ra p\\
&
\Ra g(b)\leq g(a\ra p).
\end{split}
\end{equation*}
Since $g$ is onto, it means that 
$g(a)\ra g(p)\leq g(a\ra p)$.
Conversely, 
\begin{equation*}
\begin{split}
c\leq g(a\ra p)&
\LRa g\adj(c)\leq a\ra p\\
&
\LRa g\adj(c)a\leq p\\
&
\Ra cg(a)=g(g\adj(c)a)\leq g(p)\\
&
\LRa c\leq g(a)\ra g(p),
\end{split}
\end{equation*}
i.e., $g(a\ra p)\leq g(a\ra p)$. 
Similarly we can show $g(p)\la g(a)=g(p\la a)$. 
Hence $g[P]$ is cyclic and $Q/\sim$ is simple.

The last assertion follows from 
\begin{align*}
f\adj(c)\ra f\adj(d)&
= f\adj(c\ra d),\\
f\adj(d)\la f\adj(c)&
= f\adj(d\la c).
\end{align*} 
Namely, 
\begin{equation*}
\begin{split}
a\leq f\adj(c)\ra f\adj(d)&
\LRa af\adj(c)\leq f\adj(d)\\
&
\LRa f(a)c=f(a)f(f\adj(c))=f(af\adj(c))\leq d\\
&
\LRa f(a)\leq c\ra d\\
&
\LRa a\leq f\adj(c\ra d)
\end{split}
\end{equation*}
and similarly for ``$\la$''.
\end{proof}

\begin{corol}
Let $Q$ be a quantale. Then the following sets  bijectively correspond.
\begin{enumerate}
\item
two-sided primes,
\item
cyclic sets,
\item
maximal congruences,
\item
(isomorphism classes of) simple quotients.\label{four-biject}
\end{enumerate}
\end{corol}

\begin{remark}
\begin{enumerate}
\item
A semisimple factor $Q$ is simple. Indeed, from $0\neq 1$ it follows that
$Q$ has at least one simple quotient. On the other hand, $Q$ has at most one
two-sided prime $0$. Hence the simple quotient is unique and so the 
quotient map is an isomorphism.
\item
A factor with at least one prime has a unique simple quotient. If $p$ is
the prime, then $1\x 1\x 1\not\leq p$, i.e., $1\ra p\la 1\neq 1$. Since  
$1\ra
p\la 1$ is two-sided, it must be $0$. By \ref{primes_are_closed} we have that
$0$ is a prime. It gives a simple quotient which is unique by the same argument
as in (1).
\end{enumerate}
\end{remark}

\begin{theorem}\label{semisimple}
The following conditions are equivalent:
\begin{enumerate}
\item
$Q$ is a semisimple quantale.
\item
$Q$ is strictly faithful and spatial.
\item
$Q$ is strictly faithful and $\T(Q)$ is a spatial locale.
\item
The union of cyclic sets of $Q$ is separating.
\end{enumerate}
\end{theorem}

\begin{proof}
(1)$\Ra$(2): 
A semisimple quantale is evidently spatial.
Let $a,b\in Q$ and suppose that $lar=lbr$ for all 
$l\in\mathcal{L}(Q),r\in\R (Q)$.
Then for every surjective morphism $f:Q\ra K$ onto simple quantale $K$ we
also have 
$$f(l)f(a)f(r)=f(lar)=f(lbr)=f(l)f(b)f(r)$$ 
and $\mathcal{L}(K)=f[\mathcal{L}(Q)]$, $\R(K)=f[\R (K)]$. 
Thus, by the strict faithfulness of $K$, $f(a)=f(b)$ and so $a=b$.

(2)$\Ra$(3): 
If $Q$ is spatial, $\T (Q)$ is also spatial and idempotent, 
i.e., it is a spatial locale \cite{Kr00,PaRo00}.

(3)$\Ra$(4): 
Let $a\not\geq b$. 
Then there exist $l\in\mathcal{L}(Q),r\in\R (Q)$ such that 
$lar\not\geq lbr$. 
Since $lar,lbr\in\T (Q)$ and $\T (Q)$ is spatial, there exists a two-sided 
prime $p$ with the property $lar\leq p,lbr\not\leq p$, or equivalently 
$a\leq r\ra p\la l,b\not\leq r\ra p\la l$. 
Thus we have cyclic set $\{ c\ra p\la d\mid c,d\in Q\} -\{ 1\}$ separating 
$a,b$.

(4)$\Ra$(1): 
It follows directly from \ref{simple-quotients}.
\end{proof}

\begin{remark}
In \cite{Pa97} it was shown that every simple quantale $K$ embeds into
$\Q\R(K)$ via left action on right-sided elements. Thus every surjective
morphism $f:Q\ra K$ onto simple $K$ can be extended to the quantale
representation $\mu:Q\ra\Q\R(K)$. One can easily check that this
representation is equivalent to an arbitrary representation of the form
$\mu_p:Q\ra\Q(Q/p)$ (see \cite{Kr00}) where $p$ is an element of
corresponding cyclic set $P$. More generally, $P$ itself (or its suitable
part) can be ordered into a matrix giving an equivalent representation
$\mu_P:Q\ra\Q(Q/P)$. It means that semisimplicity could be viewed as a
stricter kind of spatiality. In case of locales, every representation is
of the form $f:L\ra\Two$ and so surjective. Hence a locale is semisimple
whenever it is spatial. This result is also a consequence of the previous
theorem, because $1\wedge a\wedge 1=a$ for every $a\in L$, i.e. locales
are strictly faithful. The interpretation of a cyclic set as a (canonical)
matrix of a module will be utilized in section \ref{modules}. 
\end{remark}

\begin{example}
Let $S$ be a sup-lattice and $\Q(S)$ a quantale of its sup-lattice
endomorphisms. 
Then elements of the form
$$(\rho_x\vee\lambda_y)(z)=\begin{cases}
1 & z\not\leq y, \\
x & 0\neq z\leq y, \\
0 & z=0
\end{cases}$$
for $x\neq 1,y\neq 0$ form a cyclic separating set of $\Q(S)$. But not all
primes of $\Q(S)$ have to be of this form. For instance, if 
$S=\mathcal{M}_5$ is a five-element ``diamond'' lattice, then the identity 
satisfies $\alpha 1\beta\leq id\LRa\alpha 1\beta =0$. Since $0$ is a 
prime, $id$ must be a prime too and it is not of the above form. The same 
is hold for all endomorphisms arising from permutations of three atoms of 
$\mathcal{M}_5$.

We can construct now a proper class of non-isomorphic strong modules (see
\cite{Kr00} for definitions of used terms) over $\Q(\mathcal{M}_5)$.
Indeed, for every index set $I$ of arbitrary
cardinality we have a prime matrix with $id$ on the diagonal and $0$ otherwise.
Strong modules produced by these matrices have different cardinalities, hence
they are not isomorphic. On the other hand, every quantale has only a set of
cyclic sets. This example explains that, in some situations, the set of cyclic
sets can be a better notion of a quantale spectrum then a set of primes or a
category of strong modules.
\end{example}

\section{*-simple and *-semisimple *-quantales}

\begin{defin}
By a {\em *-simple *-quantale} is meant a *-quantale  $Q$ such that
$1\cdot 1\neq 0$ and every its quotient *-morphism is either a 
*-isomorphism or a constant *-morphism on a trivial *-quantale. 
In other words, there are exactly two *-congruences on $Q$.

By a {\em *-semisimple *-quantale} is meant a subdirect product of 
*-simple *-quant\-ales. 
That is, a *-quantale is *-semisimple if it has enough surjective 
*-morphisms onto *-simple *-quantales to separate elements.
\end{defin}

\begin{defin}
A *-quantale $Q$ is called a {\em *-factor} if $\mathcal{H}\T (Q)=\{ 
0,1\}$ and 
$0\neq 1$.
\end{defin}

\begin{defin}
Let $Q$ be a *-quantale. A set $P\se Q$ is said to be {\em *-separating} 
if $P\cup P^*$ is separating (where $P^*=\{ p^*\mid p\in P\}$).
\end{defin}

\begin{theorem}
The following conditions are equivalent:
\begin{enumerate}
\item
$Q$ is a *-simple *-quantale.
\item
$Q$ is a strictly faithful *-factor.
\item
$Q$ has a *-separating cyclic set.
\end{enumerate}
\end{theorem}

\begin{proof}
$(1)\Ra (2)$: 
If $t\in\mathcal{H}\T(Q)$, then a mapping $a\mapsto a\vee t$ gives a 
surjective 
*-morphism. Since $Q$ is *-simple, $t$ is either $0$ or $1$, i.e., $Q$ is 
*-factor. From $1\x 1\in\mathcal{H}\T(Q)$ we have $1\x 1=1$. The 
congruence defined 
in \ref{simple} is also a *-congruence, and by the same argument $Q$ is 
strictly faithful.

$(2)\Ra (3)$:
Suppose that $Q$ is strictly faithful *-factor. 
Since $1\x 1\in\mathcal{H}\T(Q)$ and $1\x 1\x 1\neq 1\x 0\x 1$, we have 
$1\x 1=1$. 
Let $u,t\in\T(Q)-\{1\}$. 
Then $uu^*,tt^*\in\mathcal{H}\T (Q)$ gives $uu^*=tt^*=0$ 
and $(u\vee t)(u^*\vee t^*)=ut^*\vee u^*t\in\mathcal{H}\T (Q)$ implies 
that 
$u\vee t=1$ or $ut^*=0$. 
Similarly, from 
$(u\vee t^*)(u^*\vee t)=ut\vee t^*u^*$ we derive
$u\vee t^*=1$ or $ut=0$. 
If both $ut=0,ut^*=0$, then $t=0$ or $t\vee t^*=1$, which implies
$u1=u(t\vee t^*)=0\Ra 1u1=0\Ra u=0$ by the strict faithfulness. 
If both $u\vee t=1,u\vee t^*=1$, then from $t\neq 1$ we have 
$$1 = (u\vee t)(u\vee t^*)
= u^2\vee ut^*\vee tu\vee tt^* 
\leq u1\vee u1\vee 1u\vee 0
\leq u,$$
a contradiction. 
If, for instance, $u\vee t=1,ut=0$, then 
\begin{align*}
u1&
=u1\x 1=u(u^*\vee t^*)1=ut^*1\leq t^*1\\
\intertext{and}
t^*1&
=t^*1\x 1=t^*(u\vee u^*)1=t^*u1\leq u1,
\end{align*} 
i.e., $u1=t^*1$. 
Now assume $ur\neq 0$ for some $r\in\R (Q)$. 
Then we have 
$$
u1r 
\leq ur
\leq u1
= u(ur\vee r^*u^*)
= u^2r
\leq u1r.
$$
Otherwise $ur=0\Ra u1r=0$. 
That is $ur=u1r$ for every $r\in\R (Q)$, and hence $lur=lu1r$ for every 
$l\in\mathcal{L}(Q),r\in\R (Q)$ gives that $u=u1$. 
Similarly it can be shown that $t^*=t^*1$, thus $u=t^*$. 
In the last case $u\vee t^*=1,ut^*=0$, by a similar trick we derive $u=t$. 
Altogether, $\T(Q)$ is a two-element chain $\{ 0,1\}$ or a four-element 
Boolean algebra $\{ 0,t,t^*,1\}$. 
Put $v=0$ or $v=t$ respectively to these two cases and 
$P=\{ r\ra v\la l \mid r\in\R(Q), l\in\mathcal{L}(Q)\} - \{1\}$. 
Then $P$ is cyclic (see \ref{simple}) and so is $P^*$. 
In case of $v=0$ the set $P$ is separating by \ref{simple}, otherwise
consider $a\not\geq b$, then $lar\not\geq lbr$ for some 
$l\in\mathcal{L}(Q),r\in\R(Q)$. 
Thus $lbr\not\leq w$ and $lar\leq w$ for some $w\in\{ v,v^*\}$. 
In other words $a\leq r\ra w\la l$, $b\not\leq r\ra w\la l$.

$(3)\Ra (1)$:
Let $t=\I P$. 
From \ref{cyclic} we have that $t=1\ra p\la 1$ for every $p\in P$. 
Further $1a1\leq p$ for some $p\in P$ implies $a\leq t$, thus 
$1\x 1\leq p\Ra 1\x 1\x 1\leq p\Ra 1\leq t$, 
resp.\ 
$1\x 1\leq p^*\Ra 1\leq t^*$, hence $1\x 1=1$. 
If $u\in\T(Q),u\neq 1$, then $u\leq p$ for some $p\in P\cup P^*$ and we 
have $u\leq 1\ra u\la 1\leq 1\ra p\la 1=t$ or $\ldots =t^*$, i.e., $t,t^*$ 
are maximal two-sided. 
Consider a *-congruence $\sim$ generated by $a\sim b$ for some 
$a\not\geq b$. 
We have some $p\in P\cup P^*$ such that $a\leq p,b\not\leq p$, suppose 
$p\in P$. 
Then $p=d\ra t\la c$ for some $c,d\in Q$. 
Put $u=1cad1,v=1cbd1$, from $a\leq d\ra t\la c$ we have $cad\leq t$ and
$u\leq 1t1\leq t$, and similarly from $b\not\leq d\ra t\la c$ we have
$cbd\not\leq t=1\ra t\la 1$, i.e., $v\not\leq t$. 
Now $u\sim v$ gives that $t=t\vee u\sim t\vee v=1$. 
Hence also $t^*\sim 1$ and 
$0 = \I P\wedge\I P^* = t\wedge t^* \geq tt^*\sim 1\x 1 = 1$, 
thus $Q$ is *-simple.
\end{proof}

We conclude a *-analogy of \ref{four-biject}.

\begin{theorem}
Let $Q$ be a *-quantale. 
Then the following sets  bijectively correspond.
\begin{enumerate}
\item
pairs of mutually *-adjoint two-sided primes,
\item
pairs of mutually *-adjoint cyclic sets,
\item
maximal *-congruences,
\item
isomorphism classes of *-simple quotients.
\end{enumerate}
\end{theorem}

Similarly as in \cite{Kr00} 5.3., we have the following statement.

\begin{theorem}
A *-quantale $Q$ is *-semisimple, iff it is semisimple as a quantale.
\end{theorem}

\begin{remark}
In \cite{Kr00}, the first author introduced a notion of a D-spatial 
*-quantale beacuse of the difference of simple and *-simple *-quantale.
We have also an opportunity to introduce a similar concept. 
Let us call a quantale {\em D-semisimple} if it is a subdirect product of
simple *-quantales. A cyclic set $P$ is called {\em *-cyclic} if $P^*=P$.
An ardent reader can find appropriate analogies of \ref{semisimple}, 
\ref{four-biject} using *-cyclic sets. One remarkable fact is that 
*-cyclic sets arises from hermitian two-sided primes.
\end{remark}

\section{Simple modules}\label{modules}

In ring theory, by simple modules are meant modules without proper
non-trivial submodules. One can easily see that any simple module $M$ has
also only two quotients---the trivial one and itself, because kernels of
the quotient mappings are just the submodules of $M$. This simple fact is
not longer true for the case of quantales, so we have to state the
following notion of simplicity.

\begin{defin}
Let $Q$ be a quantale. A $Q$-module $M$ is said to be {\em simple} if
\begin{itemize}
\item
its only submodules are $M$ and $\{ 0\}$,
\item
its only quotient modules are $M$ and $\{ 0\}$, and
\item
the action is non-trivial in the sense that $1_Q\x 1_M\neq 0_M$.
\end{itemize}
\end{defin}

Further, if $R$ is a (unital) ring and $M$ a simple faithful $R$-module,
then $R$ is a simple ring, i.e., its ideals are only $R$ and $\{ 0\}$.  
Indeed, let $x\in M,x\neq 0$ and $I\se R$ an ideal. Then $Ix=\{ ix\mid
i\in I\}$ is a submodule of $M$. Since $M$ is simple, $Ix$ is either $M$
or $\{ 0\}$ and from the faithfulness we have that $I$ is either $R$ or
$\{ 0\}$. Conversely, taking the ring itself we show that every simple
ring has a simple faithful module. Although the ``idealic'' tricks are not
disposal for the case of quantales, we are able to prove similar results.

\begin{lemma}\label{dual_faithful}
Let $Q$ be a quantale, $M$ a simple $Q$-module. Then $M\op$ is also simple 
and it is faithful whenever $M$ is so.
\end{lemma}

\begin{proof}
The first assertion is a consequence of module duality. Let $M$ be a 
faithful left $Q$-module, $a\not\geq b$ elements of $Q$. We have some 
$m\in M$ such that $am\not\geq bm$. Then $m\leq am\la a$ but $m\not\leq 
am\la b$, i.e., $M\op$ is faithful.
\end{proof}

\begin{theorem}\label{simple_module}
Let $Q$ be a quantale, $M$ a left $Q$-module. 
Then $M$ is simple iff
$Qm = \{ am \mid a\in Q\} = M$ for every $m\in M,m\neq 0$ and 
$m\la Q=\{ m\la a \mid a\in Q\} = M$ for every $m\in M,m\neq 1$.
\end{theorem} 

\begin{proof}
Let $M$ be simple, $m\neq 0$. 
Note that the set $Qm$ is a submodule of $M$. 
If $Qm=\{ 0\}$, then $\{ 0,m\}$ is a submodule with trivial action, thus 
it is proper---a contradiction. 
Hence $Qm$ is a non-trivial submodule, i.e., it is $M$. 
Since $M$ has only two quotients, $M\op$ has only two submodules and
by the same argument adapted to $M\op$ we get the second assertion.

The backward implication is evident.
\end{proof}

\begin{theorem}\label{simple_simple}
A quantale $Q$ is simple iff it has a simple faithful $Q$-module.
\end{theorem}

\begin{proof}
Let $Q$ be simple. 
We take $\R(Q)$ as a left $Q$-module. 
It is not hard to see that the assignment 
$a\sim b\LRa (\forall r\in\R(Q))(ar=br)$ defines a quantale congruence. 
From $1\x 1\neq 0$ we have that $\sim$ is not a full relation, hence 
$a\sim b\LRa a=b$. 
But it actually says that $\R(Q)$ is a faithful module. 
By \ref{simple} $Q$ is a strictly faithful factor. 
Since $1r\in\T(Q)$ for $r\in\R(Q)$, we have $1r=1$ whenever $r\neq 0$. 
Further $lr=lr1$ for $r\in\R(Q),l\in\mathcal{L}(Q)$, hence $r=r1$. 
Thus $rs=r1s=r1=r$ for $r,s\in\R(Q),s\neq 0$. 
Using \ref{simple_module} we have that $\R(Q)$ is a simple module.

Conversely, let $Q$ be a quantale and $M$ a left simple faithful 
$Q$-module.
For every $m\in M$ put $r_m=\S\{ a\in Q\mid a1_M\leq m\}$. 
From \ref{simple_module} it follows that there exists $a\in Q$ such that 
$a1_M=m$, hence $r_m1_M=m$. Since $r_m1_Q\x 1_M\leq r_m1_M=m$, $r_m$ is 
right-sided.
Similarly, for $l_m=\S\{ a\in Q\mid 0_M\la a\geq m\}$ we have 
$0_M\la a=m$ for some $a\in Q$, hence 
$0_M\la l_m =m$ and
$$
0_M\la (1_Ql_m)
= (0_M\la 1_Q)\la l_m
\geq 0_M\la l_m,
$$ 
i.e., $l_m$ is left-sided. 
If now $a\neq b$ in $Q$, then there is $m\in M$ such that $am\neq bm$ and 
thus $ar_m\neq br_m$. 
From \ref{dual_faithful} we have that $M\op$ is also faithful, thus there 
is $n\in M$ such that $n\la (ar_m)\neq n\la (br_m)$, hence $l_nar_m\neq 
l_nbr_m$, i.e., $Q$ is a strictly faithful quantale. 
Finally, let $t\in\T(Q),t\neq 0$. 
From \ref{simple_module} it follows that $1_Qm=1_M$ for every non-zero 
$m\in M$.
Since $M$ is faithful, $t1_M\neq 0$. Thus $tm\geq 1_Qt1_Qm=1_Qt1_M=1_M$ 
for every non-zero $m$, i.e., $t=1$. 
We have proved that $Q$ is factor, and hence $Q$ is a simple quantale by 
\ref{simple}.
\end{proof}

\begin{corol}
Let $Q$ be a quantale. Then (isomorphism classes of) simple $Q$-modules 
bijectively correspond to cyclic sets on $Q$.
\end{corol}

\begin{proof}
For a simple $Q$-module $M$ consider a set of elements 
$p^m_n=\S\{ a\in Q \mid an\leq m\}$ for all pairs $(m,n)\in M\times M$. 
By similar calculations as in the proof of \ref{simple_simple} we get 
$p^m_n=r_m\vee l_n$ and $r_x\ra p^m_n\la l_y=p^y_x$ whenever 
$m\neq 1,n\neq 0$, hence $P=\{ p^m_n \mid m\neq 1,n\neq 0\}$ is cyclic.

Conversely, for a cyclic set $P$ take arbitrary $p\in P$ and put 
$M=\{ a\ra p\mid a\in Q\}$. 
Then $M$ with operations $\I,\ra$ forms a left $Q$-module.
From \cite{Kr00} it follows that the dual module is 
$M\op=\{ p\la a\mid a\in Q\}$ with $\I,\la$. 
Hence the condition in \ref{simple_module} is fulfilled and $M$ is simple.

An element $n\in M$ corresponds to a subset $\{ p^m_n\mid m\in M,m\neq 1\}$. 
It is not hard to check that $a\ra p^m_n=p^m_{an}$ if $m\neq 1$, i.e., the 
module action on $M$ works as $\ra$ on $P$. Finally, the assignment 
$n\mapsto p^{\_}_n$ converts $\S$ in $M$ to $\I$ in $P$.
\end{proof}

\begin{corol}
A quantale $Q$ is semisimple iff a family of simple $Q$-modules separates 
elements, i.e., for $a\neq b$ there is simple module $M$ and $m\in M$ such 
that $am\neq bm$.
\end{corol}

\begin{remark}
If $Q$ is a *-quantale, then for every left $Q$-module $M$ there is a 
right $Q$-module $M^*$ with the same sup-lattice structure but with an 
action ``$*$'' given by $m*a=a^*m$. Using duality we have a left 
$Q$-module $(M^*)\op$. Obviously, $(M^*)\op$ is simple iff $M$ is so.
One can check that if $P$ is a cyclic set of $Q$ and defines a module 
$M$, then $P^*$ defines $(M^*)\op$. If $P$ is *-cyclic, then 
$M\cong(M^*)\op$, i.e., $M$ is a *-module (cf. \cite{Kr00}).
\end{remark}

\section{Simple quantales versus C*-algebras}

In this section, we point out which role the notion of cyclic sets can
play in theory of C*-algebras. A non-expert reader is recommended to some
monograph, e.g., \cite{KaRo97}, and a substantive paper \cite{MuPe00}.

\begin{defin}
Recall from \cite{MuPe00} that by a {\em spectrum} of a C*-algebra $A$ is 
meant a *-quantale $\Max A$ of all closed subspaces of $A$ with operations
\begin{align*}
\S X_i &=\overline{\sum X_i}, \\
XY &=\overline{\{ xy\mid x\in X,y\in Y\}}.
\end{align*}
\end{defin}

\begin{example}
Let $A=M_2(\mathbb{C})$ be an algebra of $2\times 2$ matrices over 
$\mathbb{C}$,
i.e., $A$ is a C*-algebra of all (bounded) linear operators on a  
$2$-dimensional Hilbert space. 
By easy computations we see that all proper non-zero right
ideals of $A$ are of the form
$$
R_{k,l} = 
\left\{\begin{pmatrix}
ka&kb\\
la&lb
\end{pmatrix} 
\mid a,b\in\mathbb{C}\right\}$$
for some fixed $k,l\in\mathbb{C}$. 
Similarly, proper non-zero left-sided ideals of $A$ are of the form
$$
L_{k,l} = 
\left\{\begin{pmatrix}
ka&la\\
kb&lb
\end{pmatrix} 
\mid a,b\in\mathbb{C}\right\}.$$
Thus all of these right (left) ideals are maximal and 
$R_{k,l}^*=L_{\bar{k},\bar{l}}$. 
From theory of operator algebras it follows that a given functional on a 
C*-algebra is a pure state iff its kernel is $R\vee R^*$ for some maximal 
right ideal $R$. Further, $AMA=A$ for every non-zero matrix $M$, hence 
$AXA=A$ for every non-zero (closed) subspace $X$ of $A$ and thus $A\ra 
X\la A=\{ 0\}$ for every proper subspace $X$. This guarantees that the set 
$P=\{ X\ra\{ 0\}\la Y\mid X,Y\in\Max A\}$ is cyclic. We assert that $P$ 
contains all kernels of pure functionals\footnote{See \ref{S-}.} of $A$.
Indeed, $L_{l,-k}R_{k,l}=\{ 0\}$ and we observe that $R_{k,l}$ is the 
greatest subspace with this property, i.e., $R_{k,l}=\{ 0\}\la L_{k,-l}$. 
By similar computations we get 
$$R_{k,l}\vee L_{m,n}=R_{n,-m}\ra\{ 0\}\la L_{l,-k},$$
particularly for $m=\bar{k},n=\bar{l}$ it gives a kernel of a pure state.
Further, one could easily check that every one-dimensional subspace of $A$ 
is a meet of elements of the form $R\vee L$ for some maximal right ideal 
$R$ and maximal left ideal $L$. Although $P$ is not separating set for 
$P$, from the above observation it follows that $P$ separates atoms of 
$\Max A$. So $\Max A$ is ``almost simple'',\footnote{It is mentioned in
\cite{PaRo00} that $\Max A$ contains a subspace not reachable by primes,
hence it is not spatial in the sense of \cite{PeRo97,Kr00}.} and the
algebra $A$ can be fully recovered from $\Max A$. By little bit more 
complicated calculations one can get a similar result for every 
$M_n(\mathbb{C}),n\in\mathbb{N}$.
\end{example}

Since pure states separate elements in arbitrary C*-algebra, one could
hope that its spectrum holds the above properties. Unfortunately, two
irreducible representations of a general C*-algebra may have the same
kernel but they need not be equivalent. It is the main reason why not all
pure states are contained in cyclic sets. Thus we do not know if cyclic
sets are separating for one-dimensional subspaces.

A way out could arise from the following concept.

\begin{defin}
Let $Q$ be a quantale, $S$ its multiplicative subsemigroup such that 
$0\in S$.

A set $P$ is called {\em $S$-separating} if for every 
$a,b\in S,a\not\geq b$ there is an element $p\in P$ provided that 
$a\leq p,b\not\leq p$.

A set $P$ is called {\em residually $S$-closed} if 
$a\ra p,p\la a\in P\cup\{1\}$ for every $a\in S,p\in P$.

A nonempty set $P,1\not\in P$ is called {\em $S$-cyclic} if 
$\{ a\ra p\la b \mid a,b\in S\} = P\cup\{ 1\}$ 
for every $p\in P$.
\end{defin}

\begin{example}\label{S-}
Let $A$ be a C*-algebra.  
Put 
$$S=\{\lb a\rb\mid a\in A\},$$
where $\lb a\rb$ denotes a one-dimensional subspace spanned by $a$,
$$P=\{\ker\phi(a\_b)\mid a,b\in A,\phi\text{ is a pure state}\},$$
i.e., $P$ contains subspaces of codimension 1 and 0 which are kernels of 
so called {\em pure functionals}.
From general theory of C*-algebras we know that pure states correspond to 
vector states of irreducible representations. 
More generally, pure functionals correspond to vector functionals. 
If $a,b$ are two linearly independent elements of $A$, then there exists 
an irreducible representation $\pi:A\ra\B(H)$ in which $\pi(a),\pi(b)$ are 
also independent, i.e., we can choose a basis of $H$ such that one of 
$a,b$ vanishes at some matrix element, while the second one is nonzero, 
for instance suppose 
$(\eta\mid\pi(a)\zeta)=0, (\eta\mid\pi(b)\zeta)\neq 0$ for some 
$\eta,\zeta\in H$. 
Put $\phi=(\eta\mid\pi(\_)\zeta)$ and we have that 
$\lb a\rb\leq\ker\phi,\lb b\rb\not\leq\ker\phi$, hence $P$ is $S$-separating.
One can see that 
$$\lb b\rb\ra\ker\phi\la\lb a\rb=\ker\phi(a\_b),$$ 
thus $P$ is residually $S$-closed. 
Now we will show that the set
$$P_{\phi}=\{\ker\phi(a\_b)\mid a,b\in A,\phi(a\_b)\neq 0\}$$
for some fixed pure state $\phi$ is $S$-cyclic. 
If $\phi(a\_a^*)$ is nonzero, then it is (real) multiple of some other 
equivalent pure state $\psi$ (i.e., $\ker\phi(a\_a^*)=\ker\psi$) and so 
there is some unitary $u\in A$ such that $\psi(u\_u^*)=\phi$ and thus 
$\ker\phi=\ker\phi(ua\_a^*u^*)$.
Similarly, we can find unitary $v$ such that 
$\ker\phi=\ker\phi(v^*b^*\_bv)$. 
Hence
\begin{equation*}
\begin{split}
\lb a^*u^*v\rb\ra\ker\phi(a\_b)\la\lb v^*b^*u\rb&
= \ker\phi(v^*b^*(ua\_ba^*u^*)v)\\
&
= \ker\phi(v^*b^*\_bv)\\
&
= \ker\phi.
\end{split}
\end{equation*}
We conclude that the union of $S$-cyclic sets of $\Max A$ is $S$-separating
for arbitrary C*-algebra $A$.
\end{example}

One tantalizing question left open is whether there are appropriate 
concepts of $S$-faithful, $S$-simple, $S$-semisimple quantales. 
Is there some link to simple and semisimple C*-algebras? 
We hope that an answer to this question, in either direction, will shed 
some light on the current debate about the quantale--C*-algebra  
connections.

\providecommand{\bysame}{\leavevmode\hbox to3em{\hrulefill}\thinspace}
\providecommand{\MR}{\relax\ifhmode\unskip\space\fi MR }
\providecommand{\MRhref}[2]{%
  \href{http://www.ams.org/mathscinet-getitem?mr=#1}{#2}
}
\providecommand{\href}[2]{#2}

\end{document}